\documentclass[12pt]{article}
\usepackage{amsmath,amsfonts}
\begin{document}

\baselineskip 20pt
\title{Rubel's problem on bounded analytic functions\footnote{2010 Mathematics Subject Classification: 
Primary 30H05, 30H10. Key words: bounded analytic function, Fatou point, radial limit, Rubel's ptoblem.}
}

\author{Arthur A.~Danielyan}

\maketitle 

\begin{abstract}

\noindent The paper shows that for any $G_\delta$ set $F$ of Lebesgue measure zero on the unit circle $T$ there exists
  a function $f \in H^{\infty}$ such that
the radial limits of $f$ exist at each point of $T$ and vanish precisely on $F$. 
This solves a problem proposed by Lee Rubel in 1973.

\end{abstract}

\begin{section}{Introduction.}

Let $\Delta$ and  $T$ respectively be the open unit disc and its boundary circle in the complex plane $\mathbb C$. 
As usual, we denote by $H^\infty$ the space of all bounded analytic functions in $\Delta$. 
It is well known that every $f \in H^{\infty}$ has radial limits $f(e^{i\theta})$ a. e. on $T$.
A point $e^{i\theta} \in T$ is called a Fatou point 
for $f \in H^{\infty}$ if $f(e^{i\theta})$ exists.
Below we assume that any function  $f \in H^{\infty}$ is defined also a. e. on $T$ by its radial limits 
$f(e^{i\theta}).$

The main purpose of this paper is to give  an affirmative solution to L. Rubel's Problem 5.29 published in the well-known research problem collection
of W. Hayman \cite{Haym} on the materials of ``Symposium on complex analysis" 
held in 1973 at the University of Kent, Canterbury. 
The formulation of the problem is the following.

 \vspace{0.25 cm} 
 
{\bf Problem 5.29} (See \cite{Haym}, p.168). {\it Let $F$ be a $G_\delta$ of measure zero on $T$. Then does there exist an $f \in H^{\infty}, f \neq 0,$
 such that $f=0$ on $F$ and every point of $T$  is a Fatou point of $f$?}

 \vspace{0.25 cm} 

The following minor modification of the problem asks a slightly more precise question.

\vspace{0.25 cm} 
 
{\bf Modified Problem 5.29.} Let $F$ be a $G_\delta$ of measure zero on $T$. Then does there exist an $f \in H^{\infty}$
 such that $f=0$ {\it precisely} on $F$ and every point of $T$  is a Fatou point of $f$?

 \vspace{0.25 cm} 

Problem 5.29 has remained open since it was proposed.
The following theorem completely solves both Problem 5.29 and its modification.

\vspace{0.25 cm} 
 
{\bf Theorem 1.} {\it Let $F$ be a $G_\delta$ of measure zero on $T$. Then there exists a non-vanishing $f \in H^{\infty}$ (even $\Re f >0$ on $\Delta$)
 such that $f=0$ precisely on $F$ and every point of $T$  is a Fatou point of $f$.  }

 \vspace{0.25 cm} 

Note that Theorem 1 in a sense is an extension of Fatou's following classical interpolation theorem of 1906:
 {\it If $F$ is closed and of measure zero on $T$, then there exists an element in the disc algebra which vanishes precisely on $F$}
 (see for example, \cite{Hoff}, p. 80).

 Now  assume that for some set $F \subset T$ there exists an $f \in H^{\infty}$
 such that $f=0$ precisely on $F$ and every point of $T$  is a Fatou point of $f$.
 Then $F$ is $G_\delta$ since it is the zero set on $T$ of the function $f$ which belongs 
 to the first Baire class on $T$. Also, by the classical boundary uniqueness theorem, $F$ is of measure zero on $T$. 
 Thus, Theorem 1 can be formulated also as the following ``if and only if" result.
 
 \vspace{0.25 cm} 
 
{\bf Corollary 1.} {\it  Let  $F \subset T$. There exists an $f \in H^{\infty}$ such that $f=0$ precisely 
on $F$ and every point of $T$  is a Fatou point of $f$ if and only if $F$ is 
 a $G_\delta$ of measure zero on $T$.  }

 \vspace{0.25 cm}

 As a corollary of (the proof of) Theorem 1 we also have the following description of
 the peak sets for those  elements of  $H^{\infty}$ for which
 all points of $T$ are Fatou points. 
 
\vspace{0.25 cm} 
 
{\bf Corollary 2.} {\it Let $F$ be a $G_\delta$ of measure zero on $T$.  
Then there exists a  $\lambda \in H^{\infty}$ such that: 
(a) All points of $T$ are Fatou points of $\lambda$; 
(b)  $\lambda=1$ on $F$; and (c) $|\lambda|<1$ on $\overline U \setminus F$.}

 \vspace{0.25 cm} 
 
 As above the converse implication is obvious and Corollary 2 in fact is the complete 
 description of peak sets for those  elements of  $H^{\infty}$ for which
 all point of $T$ are Fatou points. 
  
The following lemma is due to S.V. Kolesnikov (see Lemma 2 in \cite{Kol}).
 
  \vspace{0.25 cm} 
 
{\bf Lemma.} {\it  Let  $G$ be an open subset on $T$ and let $F \subset G$ be a set of measure zero on $T$.
For any $\epsilon>0$ there exists an open set $O, \ F \subset O \subset G$, and a function $g \in H^{\infty}$ 
such that: 

1) $|g(z)|<2, \ 0<\Re g(z)<1$ for $z \in \Delta$; 

2) the function $g$ has a finite radial limit $g(\zeta)$ at each point
$\zeta \in T$; 

3) at the points $\zeta \in O$ the function $g$ is analytic and $\Re g(\zeta)=1$; 

4) $|g(z)| \leq \epsilon$ on every radius 
$R_{\zeta_0}$ with end-point at $\zeta_0 \in T \setminus G$.}

 \vspace{0.25 cm} 

We use this lemma
in our proof of Theorem 1 (we repeat some relevant arguments from \cite{Kol} for the sake of completeness). 

 The main result of the paper \cite{Kol} is the following theorem of S.V. Kolesnikov, which solves the classical problem on the description of
 the sets of nonexistence of radial limits of bounded analytic functions.
 
\vspace{0.25 cm} 
 
{\bf Theorem (Kolesnikov).} {\it Let $E \subset T$.
 There exists an $f \in H^{\infty}$ such that
the radial limits of $f$ exist exactly on the set $T \setminus E$   if and only if $E$ is a $G_{\delta \sigma}$ of measure zero.}

 \vspace{0.25 cm} 
 
 The necessity part of this theorem is a well-known elementary result, while the 
 sufficiency part uses the above lemma and Carath\'{e}odory's general theorem
 on the boundary correspondence under the conformal mappings (involving the concept of a prime end).
 
In conclusion of the present paper, however, we completely eliminate  Carath\'{e}odory's theorem from
the proof of Kolesnikov's theorem. The main ingredient of this simplified proof is Kolesnikov's lemma (of course),
but  we just apply Theorem 1, which makes the presentation shorter.

\end{section}

\begin{section}{Proofs.}

{\it Proof of Theorem 1.} 
We denote by $m$ Lebesgue measure on $T$. As a $G_\delta$ of (Lebesgue) measure zero, the set $F$ is an intersection of open sets $G_k$ on $T$ such that $m(G_k)<1/2^k$, $k=1, 2, ...$. 
We assume $G_{k+1} \subset G_k$ (otherwise replace each $G_k$ by $\cap_{j=1}^{k}G_j$). 

We apply the Lemma for  $F$ and $G_k$,
and for $ \epsilon = 1/2^k$. 
Thus, we have the open sets $O_k$ on $T$,   $F \subset O_k \subset G_k$,
and the functions $g_k \in H^{\infty}$ such that for each $k$:  

(i) $|g_k(z)|<2, \ 0<\Re g_k(z)<1$ for $z \in \Delta$; 

(ii) the function $g_k$ has a finite radial limit $g_k(\zeta)$ at each point
$\zeta \in T$; 

(iii) at the points $\zeta \in O_k$ the function $g_k$ is analytic and $\Re g_k(\zeta)=1$; 

(iv) $|g_k(z)| \leq \epsilon_k$ on every radius 
$R_{\zeta_0}$ with end-point at $\zeta_0 \in T \setminus G_k$.

  \vspace{0.25 cm}  

Since by (i) each $g_k$ is bounded by $2$ and by (iv) the radial limits of $g_k$ on $T \setminus G_k$ are bounded by $1/2^k$,
by the Cauchy integral representation of the function $g_k$ we have

$$|g_k(z)| = \bigg| \frac{1}{2\pi i}\int_T\frac{g_k(\zeta)}{\zeta - z}d\zeta\bigg| \leq  \frac{1}{2\pi}\int_{T\setminus G_k} \frac{1/2^k}{|\zeta - z|}|d\zeta| + 
\frac{1}{2\pi}\int_{G_k} \frac{2}{|\zeta - z|}|d\zeta| $$
$$\leq \frac{1/2^k}{1- |z|} + 
 \frac{2}{2\pi(1 - |z|)} (1/2^k).$$
 
  \vspace{0.25 cm}  
 
This estimate clearly implies that the series $\sum_{k=1}^{\infty}g_k(z)=h(z)$ converges uniformly on
compact subsets of $\Delta$ to an analytic  function $h$ on $\Delta$ (cf. \cite{Kol}).

Since by (i) we have  $\Re g_{k}(z)>0$ for $z\in \Delta$, we also have $\Re h(z)>0$ for $z\in \Delta$. 
By (iii) we have  $\Re g_k(\zeta)=1$ on $O_k$ and thus the radial limit of $\Re h(z)$ 
is $+\infty$ at each point of $F$. 

Now let $\zeta_0 \in T \setminus F$. Then $\zeta_0 \in G_k$ only for finite many values of $k$, and by (iv), for all large enough $k$ we have 
$|g_k(z)| \leq 1/2^k$ on the radius $R_{\zeta_0}$ (with end-point at $\zeta_0$). Thus the series $\sum_{k=1}^{\infty}g_k(z)=h(z)$ converges uniformly
on the radius $R_{\zeta_0}$. Also, by (ii) each $g_k$ has a finite radial limit at $\zeta_0$ and thus $h$ has a finite radial limit at $\zeta_0$ (cf. \cite{Kol}).

The radial limit properties of the function $1+h$ are evident from above; we also note that  $\Re (1+h(z)) >1$ for  $z \in \Delta$. 
In particular, $1 + h$ has finite and nonzero radial limits everywhere on the set  $T \setminus F$. The
analytic function $f =  1/(1+h)$ is bounded by $1$ and has finite and nonzero radial limits everywhere on  $T \setminus F$. Obviously $f$ is zero free
on $\Delta$ and moreover $\Re f >0$ on $\Delta$.
Since the radial limit of $\Re h(z)$ is $+\infty$ at each point of $F$, the
radial limit of $f$ 
is zero  at each point of $F$. 

Theorem 1 is proved.

  \vspace{0.25 cm} 

{\it Proof of Corollary 2.} Let $h$ be the function from the previous proof. To complete the proof one can simply take $\lambda=h/(1+h)$.
This function clearly satisfies all the requirements of Corollary 2.

 \vspace{0.25 cm} 
  
Finally we simplify the proof of Kolesnikov's theorem by showing that it does not need to use prime ends at all.
Instead we apply an elementary (known) argument.

  \vspace{0.25 cm} 

{\it Simplified proof of Kolesnikov's theorem.}  Let
$E= \cup_{n=1}^\infty E_n,$ where each $E_n$ is a $G_\delta$ of measure zero as in Kolesnikov's theorem.
By Theorem 1 for each $E_n$ we have a function $f_n \in H^\infty$ with a positive real part, such that 
$f_n=0$ precisely on $E_n$ and every point of $T$ is a Fatou point of $f_n$. Since $\Re f_n(z)>0$
one can find an analytic function  $\log f_n(z)=\log |f_n(z)| + i \arg f_n(z)$ such that $|\arg f_n(z)| < \pi/2$ on $\Delta$.
We
have that   $\log |f_n(z)| \rightarrow -\infty$ as $z\in \Delta$ approaches radially to any point of $E_n$ and 
$\log |f_n(z)|$ has finite radial limits at each point of $T \setminus E_n$.
Obviously the radial limit of the bounded analytic function 
$$\varphi_n(z) =e^{i\log f_n(z)}=  e^{-\arg f_n(z)}[\cos(\log|f_n(z)|)+i\sin(\log|f_n(z)|)]$$ exists 
for each $\zeta \in T\setminus E_n$ and for no $\zeta  \in E_n$. Moreover, 
on the radii terminating on $E_n$ the oscillation of $\varphi_n$  is uniformly large and 
exceeds $e^{-\pi/2}$ (we use this property below).

The bounded analytic function $f(z)=\sum_{n=1}^\infty 1000^{-n} \varphi_n(z)$ has all desired properties.
At each $\zeta \in T \setminus E$ it
has a radial limit since each $\varphi_n$ does and the series converges uniformly.

 It remains to show that $f$
does not have radial limits
 on $E$.
If $\zeta_0 \in E$, let $E_m$ be the set with the smallest index $m$
such that $\zeta_0 \in E_m$. 
The partial sum $\sum_{n=1}^{m-1} 1000^{-n} \varphi_n(z)$ has a
 finite radial limit at $\zeta_0$. But at $\zeta_0$ 
  the oscillation of the radial limit of the term $1000^{-m} \varphi_m(z)$ 
  is not less than $1000^{-m} e^{-\pi/2}$ and the reminder series $\sum_{n=m+1}^\infty 1000^{-n} \varphi_n(z)$ does not exceed
 $1000^{-m}\ \frac{e^{\pi/2}}{999}$. Thus at $\zeta_0$ the oscillation of the radial limit of $f$ is larger than some positive number 
 (say, $0.5e^{-\pi/2}1000^{-m}$). The proof is over.

  \vspace{0.25 cm} 

{\bf Acknowledgment.} I would like to thank Larry Zalcman for his attention and  helpful remarks regarding this paper.

\end{section}

\begin{minipage}[t]{6.5cm}
Arthur A. Danielyan\\
Department of Mathematics and Statistics\\
University of South Florida\\
Tampa, Florida 33620\\
USA\\
{\small e-mail: adaniely@usf.edu}
\end{minipage}

\end{document}